\def\bibart#1#2#3#4#5#6#7
{\bibitem{#1}
#2, #4, #5 #6 (#3), #7}

\def\bibcoll#1#2#3#4#5#6#7#8
{\bibitem{#1}
#2, #4, in: #5, #6, #7, #3, pp. #8}

\def\bibbook#1#2#3#4#5#6
{\bibitem{#1} #2, #4, #5, #6, #3}

\def\bibdiss#1#2#3#4#5#6
{\bibitem{#1}
#2 (#3), #4, #5, #6}

\def\qed{\hfill\ \rule{2mm}{2mm} }

\def\qex{\hfill\ \vbox{\hrule\hbox{\vrule\kern4pt\vbox{\kern4pt{}
\kern4pt}\kern4pt\vrule}\hrule}}

\newcounter{casectr}

\newcounter{claimctr}

\documentclass[12pt]{article}

\usepackage{amssymb}

\newtheorem{guess}{Guess}[section]

\newtheorem{define}[guess]{Definition}
\newtheorem{prop}[guess]{Proposition}
\newtheorem{theorem}[guess]{Theorem}

\newtheorem{exam}[guess]{Example}

\usepackage[pdftitle={Automorphism Conjecture}, pdfauthor={Bernd Schroeder},
colorlinks=true, pdfstartview=FitV, linkcolor=blue, citecolor=blue, urlcolor=blue]{hyperref}

\usepackage{color}

\begin{document}

\bibliographystyle{plain}

\title{
Small Representations of Permutation Groups
by Embedding the Domain in an Ordered Set}

\author{
\small Bernd S. W. Schr\"oder \\
\small School of Mathematics and Natural Sciences\\
\small The University of Southern Mississippi\\
\small
118 College Avenue, \#5043\\
\small Hattiesburg, MS 39406\\
}

\date{\small \today}

\maketitle

\begin{abstract}

We present a representation
for permutation groups
as the automorphism group of an ordered set $U$
such that the automorphism
group's action
on a subset $T\subseteq U$
is the
permutation group itself.
For many imprimitive permutation groups, this representation
uses fewer points than the smallest representations to date.

\end{abstract}


\section{Introduction}

K\"onig's (\cite{Koenig})
classical question
when a given
group
can be represented as the
automorphism group of a (finite) graph
can be considered in an abstract and a concrete version.
The abstract version
asks if, for a
group
$G$, there is a graph
(or another combinatorial structure)
whose automorphism group (with no further structural demands made)
is isomorphic to $G$.
This question was answered affirmatively
for graphs and all finite groups
in \cite{Fruchtgraph,Fruchtdeg3,Sabidussi}.
For ordered sets and finite groups, the
question was answered in
\cite{Birkgrouprep,Thornton,BMgroupAutP}
with successively smaller underlying ordered sets
by refining Birkhoff's idea from \cite{Birkgrouprep}.
Interestingly, in \cite{Schrag},
orthomodular lattices
with a given automorphism group are constructed by using
Frucht's construction for graphs
(by way of \cite{Sabidussi}), not Birkhoff's
construction for ordered sets.
The smallest known representations
via ordered sets
have $3|G|$ elements, see \cite{BabaiDRR},
Corollary 4.3.

The concrete version of K\"onig's question
asks
whether a given
{\em permutation}
group
on $\{ 1, \ldots ,n\} $
can be the automorphism group of a graph
(or another combinatorial structure)
whose vertex
set is $\{ 1, \ldots ,n\} $.
The answer to this question is negative
for graphs (see \cite{Kagno1})
as well as for directed graphs (see \cite{Hemminggroup}).
This negative answer led to
much work to determine
{\em which}
permutation groups can be represented in this fashion,
with \cite{Barmakregular,GrechKis,GrechKisOrbit} being recent works that
also provide an overview of further research.
A permutation group is regular, iff
no nontrivial permutation has a fixed point.
In particular, the standard representation of a
group as a permutation group on itself is regular.
For regular permutation groups,
the concrete representation question was settled for
graphs in \cite{GodsilGRR}, and
for directed graphs in \cite{BabaiDRR}.
For
hypergraphs, it was addressed
and ``essentially settled" in
\cite{FoldesSinghi1,FoldesSinghi2}.

When analyzing (the number of) automorphisms of an
ordered set, it is natural to wonder if
every permutation group can occur as the action of an
automorphism group on an orbit.
Indeed, the main idea for our construction was conceived
in this context in
\cite{SchAutCon1}.
For permutation groups with a concrete hypergraph 
representation, the answer is positive by
identifying the 
hypergraph with the
ordered set that consists of the hyperedges
and of the singleton sets for the vertices, ordered by containment.
The question remains
whether, for
{\em every permutation group} $G$
on a domain $D$, there is a representation with an
ordered set $U_G $
that contains $D$ such that
$G$ is isomorphic to the automorphism group
${\rm Aut} (U_G )$
{\em and}
such that the action
${\rm Aut} (U_G )\circlearrowleft D$
of the automorphism group
${\rm Aut} (U_G )$
on the subset $D\subseteq U_G $ is $G$.
Such a representation
is given in Section \ref{repressec}.

Recall that a {\bf block}
for a permutation group
$G$ on $\{ 1, \ldots ,n \} $
is a set
$B\subseteq \{ 1, \ldots ,n \} $
such that, for all $\Phi \in G$, we have  that
$\Phi [B]\cap B\not= \emptyset $
implies $\Phi [B]= B$.
Recall that a permutation group on $\{1,\ldots ,n\} $
is called {\bf primitive}
iff the singletons and $\{1,\ldots ,n\} $
are its only blocks.
We will see in Section \ref{discusssec} that
the construction's explicit use
of an underlying block structure,
for many imprimitive permutation groups,
generates a
representation with fewer than $3|G|$
elements.

\section{The Construction}
\label{repressec}

Recall that an {\bf ordered set}
is a set $P$ equipped with a reflexive,
antisymmetric and transitive
relation $\leq $, the order relation,
and that the strict comparability
$<$ is $\leq \setminus =$.
[Note: Any inadvertently unexplained order terminology is
consistent with \cite{Schbook}.]
Recall that an {\bf antichain}
is a  subset of an ordered set such that no two
distinct elements are
comparable
and that, when an element $p$ is less than every element
of a set $X$, we write $p<X $.
Let $G\subseteq S_n $ be a permutation group.
For $A\subseteq \{ 1, \ldots , n\} $, we set
$G\cdot A=\{ \sigma [A]:\sigma \in G\} $.
Recall that the {\bf orbit} of an element $x$
is
the set $G\cdot x=\{ \sigma (x):\sigma \in G\} $, and that the
orbits partition $\{ 1, \ldots , n\} $.

\begin{define}
\label{blockpart}

Let $G$ be a subgroup of the symmetric group $S_n $
and let $\{ O_1 , \ldots , O_\ell \} $ be the partition of
$\{ 1, \ldots , n\} $ into the orbits of $G$.
For every $j\in \{ 1, \ldots , \ell \} $, let
$B_j \subseteq O_j $ be a block of $G$.
We will call the partition
${\cal B}:=\bigcup _{j=1} ^\ell G\cdot B_j $
of $\{ 1, \ldots , n\} $
a {\bf block partition} with {\bf orbit cut} $\{ B_1 , \ldots , B_\ell \} $.
For every $j\in \{ 1, \ldots , \ell \} $, we let $m_j :=|G\cdot B_j |$.
For all
$B\in {\cal B}$,
we define $G
|_{B } :=
\left\{ \theta |_{B } :\theta \in G\right\} $,

\end{define}

Note that the $\mu \in G|_{B } $ can map $B $ to any $C\in G\cdot B$.

\begin{define}
\label{nontrivblockrep}

Let $G\subseteq S_n $ be a subgroup of the symmetric group $S_n $
on
$\{ 1, \ldots , n\} =:T$,
let
${\cal B}
$
be a block partition for $G$,
and,
for all $\mu \in
\bigcup _{B\in {\cal B}} G
|_{B } $,
let
$D^\mu :={\rm dom} (\mu )\times \{ \mu \} $.
Moreover, 
let
$F:=\{ \ell _1 , u_1 , \ell _2 , u_2 ,\ldots
, \ell _n , u_n \} $.
We define
$U_{G\looparrowright {\cal B} } :=F\cup G\cup
\left( \bigcup _{\mu \in
\bigcup _{B\in {\cal B}} G
|_{B } } D^\mu \right)
\cup
T
$
to be the ordered set with the following
comparabilities.
\begin{enumerate}
\item
$\ell _1 < u_1 > \ell _2 < u_2 >\cdots
> \ell _n < u_n $, that is, $F$ is a fence with $2n $
elements.

\item
Each of the sets
$G$, $\bigcup _{\mu \in
\bigcup _{B\in {\cal B}} G
|_{B } } D^\mu $, and $T$
is an antichain.

\item
For every $\mu \in
\bigcup _{B\in {\cal B}} G
|_{B } $
and
every $j\in {\rm dom} (\mu )$, we let
$(j,\mu ) >u_j \in F$
and we let
$(j,\mu ) <\mu (j)\in T$.
\item
For
every
$\theta \in
G$,
we
let
$\theta
<\bigcup _{B\in {\cal B}} D^{\theta |_{B} }
$.
\item
Plus all comparabilities dictated by transitivity.
\end{enumerate}

\end{define}

\begin{figure}

\centerline{
\unitlength 1mm 
\linethickness{0.4pt}
\ifx\plotpoint\undefined\newsavebox{\plotpoint}\fi 
\begin{picture}(177,47)(0,0)
\put(5,30){\circle*{2}}
\put(35,30){\circle*{2}}
\put(80,15){\circle*{2}}
\put(80,5){\circle*{2}}
\put(80,45){\circle*{2}}
\put(65,30){\circle*{2}}
\put(95,30){\circle*{2}}
\put(125,30){\circle*{2}}
\put(155,30){\circle*{2}}
\put(15,30){\circle*{2}}
\put(45,30){\circle*{2}}
\put(90,15){\circle*{2}}
\put(90,5){\circle*{2}}
\put(90,45){\circle*{2}}
\put(75,30){\circle*{2}}
\put(105,30){\circle*{2}}
\put(135,30){\circle*{2}}
\put(165,30){\circle*{2}}
\put(25,30){\circle*{2}}
\put(55,30){\circle*{2}}
\put(100,15){\circle*{2}}
\put(100,5){\circle*{2}}
\put(100,45){\circle*{2}}
\put(85,30){\circle*{2}}
\put(115,30){\circle*{2}}
\put(145,30){\circle*{2}}
\put(175,30){\circle*{2}}
\put(3,28){\dashbox{1}(24,4)[cc]{}}
\put(33,28){\dashbox{1}(24,4)[cc]{}}
\put(63,28){\dashbox{1}(24,4)[cc]{}}
\put(93,28){\dashbox{1}(24,4)[cc]{}}
\put(123,28){\dashbox{1}(24,4)[cc]{}}
\put(153,28){\dashbox{1}(24,4)[cc]{}}
\put(80,15){\line(-5,1){75}}
\put(80,45){\line(-5,-1){75}}
\put(100,15){\line(5,1){75}}
\put(90,15){\line(-5,1){75}}
\put(90,45){\line(-5,-1){75}}
\put(90,15){\line(5,1){75}}
\put(100,15){\line(-5,1){75}}
\put(100,45){\line(-5,-1){75}}
\put(80,15){\line(5,1){75}}
\put(35,30){\line(3,-1){45}}
\put(35,30){\line(3,1){45}}
\put(145,30){\line(-3,-1){45}}
\put(45,30){\line(3,-1){45}}
\put(135,30){\line(-3,-1){45}}
\put(55,30){\line(3,-1){45}}
\put(125,30){\line(-3,-1){45}}
\put(65,30){\line(1,-1){15}}
\put(115,30){\line(-1,-1){15}}
\put(115,30){\line(-1,1){15}}
\put(75,30){\line(1,-1){15}}
\put(75,30){\line(1,1){15}}
\put(105,30){\line(-1,-1){15}}
\put(85,30){\line(1,-1){15}}
\put(95,30){\line(-1,-1){15}}
\multiput(45,30)(.1235955056,.0337078652){445}{\line(1,0){.1235955056}}
\multiput(90,45)(-.0786516854,-.0337078652){445}{\line(-1,0){.0786516854}}
\multiput(65,30)(.0786516854,.0337078652){445}{\line(1,0){.0786516854}}
\put(85,30){\line(-1,3){5}}
\put(95,30){\line(-1,3){5}}
\put(105,30){\line(-5,3){25}}
\multiput(125,30)(-.0786516854,.0337078652){445}{\line(-1,0){.0786516854}}
\multiput(100,45)(.0786516854,-.0337078652){445}{\line(1,0){.0786516854}}
\multiput(145,30)(-.1460674157,.0337078652){445}{\line(-1,0){.1460674157}}
\multiput(155,30)(-.1235955056,.0337078652){445}{\line(-1,0){.1235955056}}
\multiput(165,30)(-.191011236,.0337078652){445}{\line(-1,0){.191011236}}
\multiput(175,30)(-.191011236,.0337078652){445}{\line(-1,0){.191011236}}
\put(15,34){\makebox(0,0)[cb]{\footnotesize $D^{\rm id} $}}
\put(5,7){\makebox(0,0)[t]{\footnotesize ${\rm id} $}}
\put(45,34){\makebox(0,0)[cb]{\footnotesize $D^{(23)} $}}
\put(35,7){\makebox(0,0)[t]{\footnotesize ${(23)} $}}
\put(75,34){\makebox(0,0)[cb]{\footnotesize $D^{(13)} $}}
\put(65,7){\makebox(0,0)[t]{\footnotesize ${(13)} $}}
\put(105,34){\makebox(0,0)[cb]{\footnotesize $D^{(12)} $}}
\put(115,7){\makebox(0,0)[t]{\footnotesize ${(12)} $}}
\put(135,34){\makebox(0,0)[cb]{\footnotesize $D^{(123)} $}}
\put(145,7){\makebox(0,0)[t]{\footnotesize ${(123)} $}}
\put(165,34){\makebox(0,0)[cb]{\footnotesize $D^{(132)} $}}
\put(175,7){\makebox(0,0)[t]{\footnotesize ${(132)} $}}
\put(80,5){\line(0,1){10}}
\put(80,15){\line(1,-1){10}}
\put(90,5){\line(0,1){10}}
\put(90,15){\line(1,-1){10}}
\put(100,5){\line(0,1){10}}
\put(90,47){\makebox(0,0)[cb]{\footnotesize $T$}}
\put(90,3){\makebox(0,0)[ct]{\footnotesize $F$}}
\put(5,30){\line(0,-1){20}}
\put(175,30){\line(0,-1){20}}
\put(35,30){\line(0,-1){20}}
\put(145,30){\line(0,-1){20}}
\put(65,30){\line(0,-1){20}}
\put(115,30){\line(0,-1){20}}
\put(5,10){\line(1,2){10}}
\put(175,10){\line(-1,2){10}}
\put(35,10){\line(1,2){10}}
\put(145,10){\line(-1,2){10}}
\put(65,10){\line(1,2){10}}
\put(115,10){\line(-1,2){10}}
\put(25,30){\line(-1,-1){20}}
\put(155,30){\line(1,-1){20}}
\put(55,30){\line(-1,-1){20}}
\put(125,30){\line(1,-1){20}}
\put(85,30){\line(-1,-1){20}}
\put(95,30){\line(1,-1){20}}
\put(5,10){\circle*{2}}
\put(175,10){\circle*{2}}
\put(35,10){\circle*{2}}
\put(145,10){\circle*{2}}
\put(65,10){\circle*{2}}
\put(115,10){\circle*{2}}
\end{picture}
}

\caption{
Hasse diagram of the ordered set
$U_{
S_3 \looparrowright \{ \{ 1,2,3\} \} } $.
}
\label{S3_repres}

\end{figure}
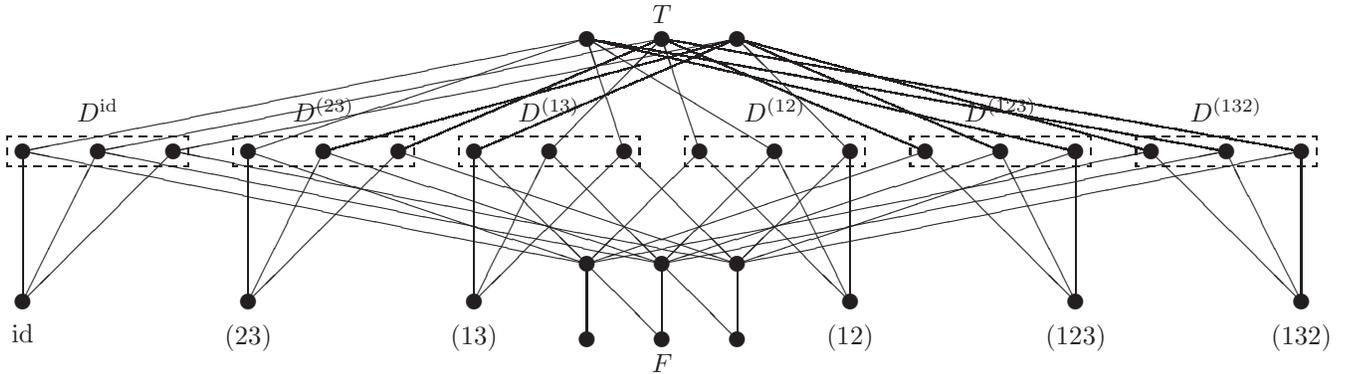

\begin{exam}

{\rm
The simplest visualization for
ordered sets
$U_{G\looparrowright {\cal B} } $
is obtained for the trivial block partition
${\cal B}=\{ \{ 1, \ldots ,n \} \} $.
In this case,
for
every
$\theta \in G$, we have that
$D^\theta  =D^{\theta |_{\{ 1, \ldots , n\} } } $
is an antichain with $n$ elements,
and
$\theta $ is simply a lower bound of
$D^\theta $.

Figure \ref{S3_repres} shows the resulting ordered set for the
symmetric group $S_3 $.
The sets $D^\theta $ are
indicated with dashed boxes.
As the proof of Theorem \ref{repwithblocks} below
will show, every permutation $\sigma $ of $T$
is induced by an automorphism that fixes every point in $F$:
Every $\theta  $ is mapped to $\sigma \theta $, and 
$D^\theta $ is mapped
$D^{\sigma \theta } $, with $(j,\theta )\mapsto (j,\sigma \theta )$.
Moreover, there are no further automorphisms.

Ordered sets
$U_{
S_n \looparrowright \{ \{ 1,\ldots ,n\} \} } $
are obtained similarly.
For subgroups
$G$ of $S_n $, ordered sets $U_{G\looparrowright \{ \{ 1,\ldots ,n\} \} }$
are obtained by, for all $\theta \not\in G$, deleting
$\theta $
and $D^\theta $ from
$U_{
S_n \looparrowright \{ \{ 1,\ldots ,n\} \} } $.
\qex
}

\end{exam}

Recall that $x$ is a {\bf lower cover}
of $y$ and $y$ is an {\bf upper cover}
of $x$ iff $x<y$ and there is no
$z\in P\setminus \{ x,y\} $
such that $x<z<y$.
Recall that a {\bf minimal element}
is an element $b$ such that there is no
$x$ with $x<b$, that
a {\bf maximal element}
is an element $u$ such that there is no
$x$ with $x>u$,
that a {\bf chain}
is an ordered set such that any
two elements are comparable, that the {\bf length}
of a chain is one less than its number of elements,
and that the {\bf rank}
of an element  $x$ is the length of a longest chain
from $x$ to a minimal element.

\begin{theorem}
\label{repwithblocks}

Let $G$ be a subgroup of $S_n $ and let
${\cal B}
$
be a block partition for $G$.
The automorphism group
${\rm Aut} \left(
U_{
G\looparrowright {\cal B} }
\right) $
is isomorphic to $G$
via the restriction $\Psi \mapsto \Psi |_T $.

\end{theorem}

{\bf Proof.}
Clearly, because
$(\Phi \Psi )|_T
=\Phi |_T \Psi |_T $,
the function $\Psi \mapsto \Psi |_T $ is
a group homomorphism.

{\em
Claim.
Every
$\Psi \in {\rm Aut} \left(
U_{
G\looparrowright {\cal B} }
\right) $
is uniquely determined by $\sigma :=\Psi |_T
$, which is an element of $G$.
}

Let $\Psi \in {\rm Aut} \left(
U_{
G\looparrowright {\cal B} }
\right) $
and let $\sigma :=\Psi |_T $.
Because automorphisms
map maximal elements to maximal elements,
$\Psi $ maps $T$ to itself.
Because automorphisms preserve the lengths of chains
and map minimal elements to minimal
elements,
$\Psi $ maps each of the sets
$\{ \ell _1 ,\ell _2 ,\ldots ,\ell _n \} $
and
$G$
to itself.
Because automorphisms preserve
the rank,
$\Psi $ maps
each of $\{  u_1 ,u_2 ,\ldots ,u_n \} $
and
$\bigcup _{\mu \in \bigcup _{B\in {\cal B} } G
|_{B } } D^\mu $
to itself.
Hence $\Psi [F]=F$, and, because $F$
has exactly one automorphism,
$\Psi |_F ={\rm id}_F $.

For the following, let
$\mu \in \bigcup _{B\in {\cal B} } G
|_{B} $.
Clearly, there is an $A\in {\cal B} $
such that $\mu \in G|_{A } $, that is, $A$ is the domain of $\mu $.

First let $j\in {\rm dom} (\mu )=A $.
Then there are a $\nu \in \bigcup _{B\in {\cal B} } G
|_{B } $ and an $i\in {\rm dom} (\nu )$ such that
$\Psi (j,\mu )=(i,\nu )\in D^\nu $.
Now
$u_j
< (j,\mu )
$
implies
$u_j =
\Psi ( u_j )
< \Psi \left( j, \mu \right)
=
(i,\nu )$.
Because
$(j,\nu )$ is the unique
upper cover of $u_j $ in $D^\nu $,
we
obtain $i=j$, that is
$\Psi \left( j, \mu \right)
=(j,\nu )$.
Because $j\in {\rm dom} (\nu )$
and
${\cal B} $ is a partition,
by definition, we obtain that
${\rm dom} (\nu )=A $.

Next, let
$(a,\mu ),(b,\mu )\in D^\mu $
and consider
$\Psi (a,\mu )=:(a , \nu _a )$
and
$\Psi (b,\mu )=:(b , \nu _b )$.
By definition,
$(a,\mu )$ and $(b,\mu )$ have the same
lower covers in $G$.
Hence
$\Psi (a,\mu )=(a , \nu _a )$
and
$\Psi (b,\mu )=(b , \nu _b )$
have the same
lower covers in $G$.
Thus, for all $\theta \in G$ with $\theta <\nu _a ,\nu _b $, we have that
$\nu _a =\theta |_{A } =\nu _b $.
Because $(a,\mu ),(b,\mu )\in D^\mu $
were arbitrary,
we conclude that
there is a
$\nu \in G
|_{A } $ such that
$\Psi \left[ D^\mu \right] =D^\nu $.

For every
$j\in A $,
we have
that $u_j <(j,\mu )<\mu (j)$
implies
$u_j
=\Psi (u_j )
<
\Psi  (j, \mu )
=
(j, \nu )
=\Psi  (j, \mu )
<
\Psi \left( \mu (j) \right) =
\sigma \left( \mu (j) \right) $.
By definition, we must have $\nu =\sigma \mu $.
Therefore
$\Psi \left[ D^\mu \right] =D^{\sigma \mu }$, and, for all
$(j,\mu )\in A \times \{ \mu \} $, we have
$\Psi (j,\mu )=(j,\sigma \mu )$.

Because the above holds for all $\mu \in \bigcup _{B\in {\cal B}}  G
|_{B } $, for every
$\theta \in G$,
we have that
$\Psi (\theta )$ is the lower cover
of $\bigcup _{B\in {\cal B}} D^{\sigma \theta |_{B } } $.
Consequently
$\Psi (\theta )
=
\sigma \theta
$.
On one hand, this concludes the proof
that $\sigma =\Psi |_T $
determines $\Psi (x)$ for every
$x\in
U_{
G\looparrowright {\cal B} }
$.
On the other hand,
using any $\theta \in G$, we have that $\sigma \theta =\Psi (\theta )\in  G$
implies $\sigma =\Psi (\theta )\theta ^{-1} \in G$.
This completes the proof of the {\em Claim}.
In particular,
$\Psi \mapsto \Psi |_T $ is
injective and it maps into $G$.

Finally, for surjectivity, let
$\sigma \in G $.
We define
$\Phi _\sigma $
by
$\Phi _\sigma |_F :={\rm id}|_F $,
by, for every $j\in T$,
setting
$\Phi _\sigma (j):=\sigma (j)$,
for every $\mu \in \bigcup _{B\in {\cal B} } G
|_{B } $
and every $j\in {\rm dom} (\mu) $,
setting
$\Phi _\sigma \left( j, \mu \right)
:=
(j, \sigma \mu )$,
and,
for every
$\theta \in G$,
setting
$
\Phi _\sigma (\theta )
=
\sigma \theta
$.
Clearly, $\Phi _\sigma $
is injective,
hence surjective, $\Phi _\sigma $
preserves
order on $F$,
and
every comparability
$\theta <D^{\theta |_{B } } $
is mapped to
$
\Phi _\sigma (\theta )
=
\sigma \theta
<D^{\sigma \theta |_{B } }
=\Phi _\sigma \left[ D^{\theta |_{B } } \right]$.
Finally,
for every $\mu \in \bigcup _{B\in {\cal B} } G
|_{B} $
and every $j\in {\rm dom} (\mu) $,
we have that $u_j< (j, \mu )<\mu  (j) $
implies
$\Phi _\sigma (u_j )
=u_j
< (j, \sigma \mu )
=
\Phi _\sigma \left( j, \mu \right)
=
(j, \sigma \mu )
<\sigma \left( \mu (j) \right)
=
\Phi _\sigma \left( \mu (j) \right)
$.
Hence $\Phi _\sigma \in {\rm Aut} (U_{G\looparrowright {\cal B} } )$
and
$\Psi \mapsto \Psi |_T $ is
surjective.
\qed

\vspace{.1in}

Now that we have a representation, we turn to
its size.

\begin{define}

Let $G$ be a subgroup of $S_n $ and let
${\cal B}
$
be a block partition for $G$.
For every
$B\in {\cal B} $, we define
$G\circlearrowleft B :=
\{ \sigma |_{B } :\sigma \in G, \sigma [B ]=B \} $.
Moreover, we define
$G\circlearrowleft {\cal B}$
to be the group of permutations on ${\cal B}$
induced by $G$ by, for every $\sigma \in G$,
mapping each block to its image
under $\sigma \in G$.

\end{define}

\begin{prop}
\label{size}

Let $G$ be a subgroup of $S_n $, let
${\cal B}
$
be a block partition for $G$,
and let
$\left\{ B_{1 } , \ldots , B_{\ell } \right\} \subseteq {\cal B}$
be an orbit cut.
Then
$$\displaystyle{
\left|
U_{
G\looparrowright {\cal B}
}
\right|
=
|G|
+
3n
+
\sum _{i=1}
^\ell
|B_{i } |\cdot
m_i ^2 \cdot |G\circlearrowleft B_{i } |
.} $$
In case
$G\circlearrowleft {\cal B}$ is transitive,
with
$N:=\left\{
\sigma \in G: \sigma [B_i ]=B_i , i=1,\ldots ,m
\right\} $,
we obtain
$\displaystyle{
\left|
U_{
G\looparrowright {\cal B}
}
\right|
=
|G|
\left(
1+
{
3n
\over
|G|
}
+
{|G\circlearrowleft B_1 |
\over |N|}
{
|{\cal B} |n
\over
|G\circlearrowleft {\cal B} |
}
\right)
.} $
In case $G$ is transitive, we obtain
$\displaystyle{
\left|
U_{
G\looparrowright {\cal B}
}
\right|
\leq
|G|
\left(
1+
{|G\circlearrowleft B_1 |
\over |N|}
\left(
n +3
\right)
\right)
.}$

\end{prop}

{\bf Proof.}
We have
$|T|=n$ and $|F|=2n $.

Now let $i\in \{ 1, \ldots , \ell \} $.
Then, because every $\mu \in G|_{B_i } $ maps onto exactly one
$B\in G\cdot B_i $, we have
$\left| G
|_{B_{i} } \right|
=
m_i \cdot \left| G\circlearrowleft B_{i} \right|
$.
Moreover, for every $B\in G\cdot B_{i} $,
and every $\mu \in G|_B $, we have
$\left| D^\mu \right| =|B_{i} |$.
Therefore, because every
$B\in {\cal B}$ is in
exactly one $G\cdot B_{i } $, we obtain
$\left|
\bigcup _{\mu \in \bigcup _{B\in {\cal B} } G
|_{B} } D^\mu
\right|
=
\sum _{i=1}
^\ell
|B_{i } |\cdot
m_i \cdot
\left| G
|_{B_{i } } \right|
=
\sum _{i=1}
^\ell
|B_{i } |\cdot
m_i \cdot
m_i \cdot |G\circlearrowleft B_{i } |
$.
Hence
$
\left|
U_{
G\looparrowright {\cal B}
}
\right|
=
|G|
+
3n
+
\sum _{i=1}
^\ell
|B_{i } |\cdot
m_i ^2 \cdot |G\circlearrowleft B_{i } |
$.

In case $G\circlearrowleft {\cal B}$ is transitive,
we have
$\ell =1 $, $m_1 =|{\cal B} |$
and
$m_1 |B_1 |=n$.
Moreover,
$N$ is a normal subgroup of
$G$ and $G/N$ is isomorphic to
$G\circlearrowleft {\cal B}$,
so $|G|=|N|\cdot |G\circlearrowleft {\cal B}|$.
Now
\begin{eqnarray*}
\left|
U_{
G\looparrowright {\cal B}
}
\right|
& = &
|G|
+
3n
+
|B_{1 } |\cdot
m_1 \cdot m_1 \cdot |G\circlearrowleft B_{1 } |
\\
& = &
|G|
+
3n
+
|{\cal B} | n
|G\circlearrowleft B_{1 } |
\\
& = &
|G|
+
|G|
{3n+|{\cal B} |n |G\circlearrowleft B_1 |
\over |G|}
\\
& = &
|G|
\left(
1+
{
3n
\over
|G|
}
+
{|G\circlearrowleft B_1 |
\over |N|}
{
|{\cal B} |n
\over
|G\circlearrowleft {\cal B} |
}
\right)
.
\end{eqnarray*}

Finally, in case
$G$ is transitive, we have
$|G\circlearrowleft {\cal B} |\geq |{\cal B}|$,
$|G|\geq |N|\cdot |{\cal B} |$,
and
$|G\circlearrowleft B_1 |\geq |B_1| ={n\over |{\cal B} |} $.
Hence we obtain
$$
{
3n
\over
|G|
}
+
{|G\circlearrowleft B_1 |
\over |N|}
{
|{\cal B} |n
\over
|G\circlearrowleft {\cal B} |
}
\leq
{
3n
\over
|N||{\cal B} |
}
+
{|G\circlearrowleft B_1 |
\over |N|}
{
|{\cal B} |n
\over
|{\cal B} |
}
\leq
{|G\circlearrowleft B_1 |
\over |N|}
\left(
n +3
\right)
$$
and the final inequality follows.
\qed

\section{Discussion}
\label{discusssec}

To see that, for
many
imprimitive permutation groups,
Theorem \ref{repwithblocks}
gives a very compact representation,
consider the
terms in Proposition \ref{size}.
For many permutation groups, we have
$|G|\geq 3n$, and
we often have
$|G|\gg 3n$.
Moreover,
${
|G\circlearrowleft B_1 |
\over
|N|
}
\leq
1$ and, for imprimitive $G$, we often have ${
|G\circlearrowleft B_1 |
\over
|N|
}
\ll
1$.
Finally,
for transitive actions
$G\circlearrowleft {\cal B} $, we have
$|G\circlearrowleft {\cal B} |\geq |{\cal B} |$
and often
$|G\circlearrowleft {\cal B} |\gg |{\cal B} |$.
Consequently,
a demand
$
{
3n
\over
|G|
}
+
{|G\circlearrowleft B_1 |
\over |N|}
{
|{\cal B} |n
\over
|G\circlearrowleft {\cal B} |
}
<1$
is a relatively mild condition, which, via
Theorem \ref{repwithblocks}
and
Proposition \ref{size},
yields a representation
with at most
$2|G|$ points.
For transitive
imprimitive groups,
with sufficiently many independent actions on the blocks
$B\in {\cal B} $, we obtain
${|G\circlearrowleft B_1 |
\over |N|}
\left(
n +3
\right)
<1$, and again
we have
a representation with
fewer than $2|G|$ points.

In fact, it is an easy task to
construct sequences
$\{ G_k \} _{k=1} ^\infty $
of transitive imprimitive permutation groups
$G_k $ whose
order goes to infinity and such that
there are nontrivial block partitions
${\cal B}_k $
such that
the quotient
${
\left|
U_{
G_k \looparrowright {\cal B}_k
}
\right|
\over |G_k |}
$
converges to $1$.
For example, if $G_k $
is the group on $\left\{ 1, \ldots , k^2 \right\} $, generated
by the cycles
$(1+\ell k, 2+\ell k, \ldots , k+\ell k)$,
$\ell =0, \ldots , k-1$,
plus any transitive action on the block partition
${\cal B} _k :=\left\{ \left\{
1+\ell k, 2+\ell k, \ldots , k+\ell k \right\} :
\ell =0, \ldots , k-1\right\} $, then
$1+
{|G\circlearrowleft B_1 |
\over |N|}
\left(
n +3
\right)
\leq
1+
{
k
\over
k^k
}
\left(
k^2 +3
\right)
\to 1$
(indeed quite rapidly)
as $k\to \infty $.
This simple example
indicates how
sufficiently many independent actions on the blocks lead to
a rather compact representation:
Indeed, with the partition above and cyclic action on $B_1 $, 
$|N|\geq k^4 $ would still suffice.

Finally, via the trivial block partition
$\{ \{ 1,\ldots ,n\} \} $,
Theorem \ref{repwithblocks}
allows
for a ``semi-concrete" representation of
{\em any} permutation
group
$G\subseteq S_n $
via
${\rm Aut} \left(
U_{
G\looparrowright \{ \{ 1,\ldots ,n\} \} }
\right) $.
However,
by Proposition \ref{size},
$
\left|
U_{
G\looparrowright \{ \{ 1,\ldots ,n\} \}
}
\right|
=
(n+1)|G|
+
3n
>3|G|$.
It must also be noted that, when $n=|G|$, which is the
case for the standard, though combinatorially voluminous,
representation of an abstract group as a permutation group
on itself,
$U_{
G\looparrowright \{ \{ 1,\ldots ,n\} \}
}
$ makes for a rather large representation.
On the positive side,
$U_{
G\looparrowright \{ \{ 1,\ldots ,n\} \}
}$
can easily be turned into a lattice
with
$(n+1)|G|
+
3n
+
3
$
elements and
the same property
as described in
Theorem \ref{repwithblocks}
by
adding a top element ${\bf 1}$,
a bottom element ${\bf 0}$,
and a ``center element" $c$ such that
$F<c<T$,
and, for all $\theta \in G$, we have
$\theta <c<T$.

\vspace{.1in}

{\bf Acknowledgment.}
The author very much thanks Frank a Campo for
reading an earlier version of this manuscript and for
identifying a way to eliminate an
unnecessary technical condition in
the construction.

\section{Declarations}

\begin{itemize}
\item
{\bf Ethical Approval.}
Not applicable, because no human or animal studies were involved.

\item
{\bf Competing Interests.}
The author has no relevant financial or non-financial
interests to disclose.

\item
{\bf Authors' Contributions.}
This is a single author paper with all
work done by the author.
Contributions by non-authors are acknowledged above.

\item
{\bf Funding.}
No funds, grants, or other support was received.

\item
{\bf Availability of Data and Materials.}
Data sharing is not applicable, as no datasets
were generated or analysed herein.

\end{itemize}

\end{document}